\documentclass[12pt, reqno]{amsart}
\usepackage{amsfonts}
\usepackage{amssymb}
\usepackage{amsmath}

\newtheorem{theorem}{Theorem}
\theoremstyle{plain}

\newtheorem{corollary}[theorem]{Corollary}

\begin{document}
\title[] {Some identities on the $q$-Bernstein polynomials, $q$-Stirling numbers and $q$-Bernoulli numbers}
\author{Taekyun Kim}
\address{Taekyun Kim. Division of General Education-Mathematics \\
Kwangwoon University, Seoul 139-701, Republic of Korea  \\}
\email{tkkim@kw.ac.kr}
\author{Jongsung Choi}
\address{Jongsung Choi. Division of General Education-Mathematics \\
Kwangwoon University, Seoul 139-701, Republic of Korea  \\}
\email{jeschoi@kw.ac.kr}
\author{Young-Hee Kim}
\address{Young-Hee Kim. Division of General Education-Mathematics\\
Kwangwoon University, Seoul 139-701, Republic of Korea  \\}
\email{yhkim@kw.ac.kr}
\maketitle

{\footnotesize {\bf Abstract} \hspace{1mm} {
In this paper, we consider the $q$-Bernstein polynomials on $\Bbb Z_p$ and
investigate some interesting properties of $q$-Bernstein polynomials related to
$q$-Stirling numbers and Carlitz's type $q$-Bernoulli numbers.
}

\section{Introduction}

Let $C[0, 1]$ denote the set of continuous functions on $[0, 1]$. Then Berstein operator for
$f \in C[0, 1]$ is defined as
$$ \mathbb{B}_n(f)(x)= \sum_{k=0}^{n} f(\frac{k}{n})\binom{n}{k} x^k (1-x)^{n-k}=\sum_{k=0}^{n} f(\frac{k}{n})B_{k,n}(x),$$
for $k, \, n \in \Bbb Z_+$, where $B_{k, n} (x)=\binom{n}{k}x^k (1-x)^{n-k}$ is called the Bernstein polynomial of degree $n$
(see [1, 9, 10]).

In [9], Phillips introduced the $q$-extension of Bernstein polynomials and
Kim-Jang-Yi proposed the modified $q$-Bernstein polynomials of degree $n$,
which are different $q$-Berstein polynomials of Phillips (see [1]).

Let $q$ be regarded as either a complex number $q\in\Bbb C$ or a $p$-adic number
$q\in\Bbb C_p$. If $q\in\Bbb C$, then we always assume that $|q|< 1$.
If $q\in\Bbb C_p$, we normally assume that $|1-q|_p < p^{- \frac{1}{p-1}}$,
which yields the relation $q^x = \exp(x \log q)$ for $|x|_p \le 1$ (see [1-11]).

Here, the symbol $| \cdot |_p$ stands for the $p$-adic absolute on $\Bbb C_p$
with $|p|_p = 1/p$. The $q$-bosonic natural numbers are defined by
$[n]_q=\frac{1-q^n}{1-q} =1+q+ \cdots + q^{n-1} \, (n \in \Bbb N)$, and the $q$-factorial is defined by
$[n]_q!=[n]_q[n-1]_q\cdots[2]_q[1]_q$.  In this paper we use the notation
for the Gaussian binomial coefficients in the form of
$${\binom{n}{k}}_q=\frac{[n]_q!}{[n-k]_q![k]_q!}=\frac{[n]_q[n-1]_q\cdots[n-k+1]_q}{[k]_q!}, \quad (\text{see [3, 8]}).$$
Note that $\underset{q \rightarrow 1}{\lim} {\binom{n}{k}}_q ={\binom{n}{k}} = \frac{n(n-1)\cdots(n-k+1)}{n!}.$

Let $p$ be a fixed prime number. Throughout this paper, the symbol
$\Bbb Z_p$, $\Bbb Q_p$, $\Bbb C$ and $\Bbb C_p$ denote the ring of $p$-adic integers, the field of
$p$-adic rational numbers, the complex number field, and the completion of algebraic closure
of $\Bbb Q_p$, respectively. Let $\mathbb{N}$ be the set of natural
numbers and $\Bbb Z_+ = \Bbb N \cup \{ 0 \}$.

Let $UD(\Bbb Z_p)$ be the space of uniformly differentiable function
on $\Bbb Z_p$. For $f \in UD(\Bbb Z_p)$, the $p$-adic $q$-integral on $\Bbb Z_p$ is defined by
$$I_q (f)=\int_{\Bbb Z_p
} f(x) d\mu_q (x) = \lim_{N \rightarrow \infty} \frac{1}{[p^N]_q} \sum_{x=0}^{p^N
-1} f(x)q^x, \quad (\text{see [4, 5, 6]}).
$$
Carlitz' $q$-Bernoulli number $\beta_{k,q}$ can be defined respectively by
$\beta_{0,q}=1$ and by the rule that $q(q \beta +1)^k -  \beta_{k,q}$ is equal to $1$ if $k=1$ and to $0$ if $k>1$
with the usual convention of replacing $\beta^i$ by $\beta_{i, q}$ (see [5]). As was shown in [5],
Carlitz's $q$-Bernoulli numbers can be represented by $p$-adic $q$-integral on $\Bbb Z_p$ as follows :

\begin{eqnarray}\int_{\Bbb Z_p
} [x]_q^n d\mu_q (x) = \lim_{N \rightarrow \infty} \frac{1}{[p^N]_q} \sum_{x=0}^{p^N
-1} [x]_q^n q^x = \beta_{n, q}, \quad n \in \Bbb Z_+.
\end{eqnarray}

The $k$-th order factorial of the $q$-number $[x]_q$, which is defined by
$$[x]_{k,q} =[x]_q [x-1]_q \cdots [x-k+1]_q = \frac{(1-q^x)(1-q^{x-1})\cdots(1-q^{x-k+1})}{(1-q)^k}, $$
is called the $q$-factorial of $x$ of order $k$ (see [3]). Thus, we note that $\binom{x}{k}_q = \frac{[x]_{k,q}}{[k]_q !}.$

In this paper, we consider $q$-Bernstein polynomials on $\Bbb Z_p$ and we investigate some interesting properties of
$q$-Bernstein polynomials related $q$-Stirling numbers and Carlitz's $q$-Bernoulli numbers.

\section{$q$-Bernstein polynomials related to $q$-Stirling numbers and $q$-Bernoulli numbers}

In this section, we assume that $q \in \Bbb C_p$ with $|1-q|_p < p^{- \frac{1}{p-1}}$.
For $f \in UD(\Bbb Z_p)$, we consider $q$-Bernstein type operator on $\Bbb Z_p$ as  follows :
\begin{eqnarray} \mathbb{B}_{n,q}(f)(x)&=& \sum_{k=0}^{n} f(\frac{k}{n})\binom{n}{k} [x]_q^k [1-x]_q^{n-k}\\
&=&\sum_{k=0}^{n} f(\frac{k}{n})B_{k,n}(x, q),\notag \end{eqnarray}
for $k, \, n \in \Bbb Z_+$, where $B_{k, n} (x, q)=\binom{n}{k}[x]_q^k [1-x]_q^{n-k}$ is called q-Bernstein type polynomials of degree $n$
(see [1]).

Let $(Eh)(x)=h(x+1)$ be the shift operator. Consider the $q$-difference operator as follows :
\begin{eqnarray}
\Delta_q^n = \prod_{i=1}^n (E- q^{i-1}I),
\end{eqnarray}
where $(Ih)(x)=h(x)$. From (3), we note that
\begin{eqnarray}
f(x)=\sum_{n \ge 0} \binom{x}{n}_q \Delta_q^n f(0), \quad (\text{see [3]}),
\end{eqnarray}
where
\begin{eqnarray}
\Delta_q^n f(0)=\sum_{k=0}^{n}  \binom{n}{k}_q (-1)^k q^{\binom{k}{2}} f(n-k).
\end{eqnarray}
The $q$-Stiring number of the first kind is defined by
\begin{eqnarray}
\prod_{k=1}^n (1+[k]_q z) = \sum_{k=0}^{n} S_1 (n, k : q)z^k,
\end{eqnarray}
and the $q$-Stiring number of the second kind is also defined by
\begin{eqnarray}
\prod_{k=1}^n (\frac{1}{1+[k]_q z}) = \sum_{k=0}^{n} S_2 (n, k : q)z^k, \quad (\text{see [3]}).
\end{eqnarray}
By (3), (4), (5), (6) and (7), we see that
\begin{eqnarray}
S_2 (n, k : q)z^k = \frac{q^{- \binom{k}{2}}}{[k]_q !} \sum_{j=0}^{k}(-1)^j q^{\binom{j}{2}}\binom{k}{j}_q [k-j]_q^n
= \frac{q^{- \binom{k}{2}}}{[k]_q !} \Delta_q^k 0^n, \quad (\text{see [1]}).
\end{eqnarray}

For $q \in \Bbb C_p$ with $|1-q|_p < p^{- \frac{1}{p-1}}$, we have that for $k, n \in \Bbb{Z_+}$,
\begin{eqnarray*}
F_q^{(k)}(t,x)&=& \frac{t^k e^{[1-x]_q t} [x]_q^k}{k!}=[x]_q^k \sum_{n=0}^{\infty}
\binom{n+k}{k} [1-x]_q^n \frac{t^{n+k}}{(n+k)!} \\
&=&\sum_{n=k}^{\infty}
\binom{n}{k}[x]_q^k [1-x]_q^{n-k} \frac{t^n}{n!}\\
&=&\sum_{n=0}^{\infty} B_{k,n}(x,q) \frac{t^n}{n!}, \quad \text{(see [1])}.
\end{eqnarray*}
 Thus, we note that
$\frac{t^k e^{[1-x]_q t} [x]_q^k}{k!}$ is the generating function of $q$-Berstein polynomials (see [1]).
It is easy to show that
\begin{eqnarray}
[1-x]_q^{n-k} = \sum_{m=0}^{\infty} \sum_{l=0}^{n-k} \binom{l+m-1}{m} \binom{n-k}{l} (-1)^{l+m} q^l [x]_q^{l+m}(q-1)^m.
\end{eqnarray}
By (1), (2) and (9), we obtain the following theorem.

\begin{theorem}
For $k, n \in \Bbb{Z_+}$ with $n \ge k$, we have
\begin{eqnarray*}
& & \int_{\Bbb Z_p} B_{k, n}(x, q)d \mu_q (x) \\ & & \qquad=\sum_{m=0}^{\infty} \sum_{l=0}^{n-k}\binom{l+m-1}{m} \binom{n-k}{l}
(-1)^{l+m} q^l (q-1)^m \beta_{l+m+k, q},
\end{eqnarray*}
where $\beta_{n,q}$ are the $n$-th Carlitz $q$-Bernoulli numbers.
\end{theorem}

\medskip

In [3], it is known that
\begin{eqnarray}
[x]_q^n = \sum_{k=0}^{n} q^{\binom{k}{2}} \binom{x}{k}_q [k]_q !S_2 (n, k : q),
\end{eqnarray}
and
\begin{eqnarray}
\frac{\sum_{k=i-1}^{n} \frac{\binom{k}{i}}{\binom{n}{i}}B_{k,n} (x,q)}{([x]_q +[1-x]_q)^{n-i}}
=[x]_q^i, \quad \text{for} \,\, i \in \Bbb N, \quad (\text{see [1]}).
\end{eqnarray}
By (2) and (11), we see that for $i \in \Bbb N$,
\begin{eqnarray}
&&[x]_q^i  = \sum_{m=0}^{\infty} \sum_{k=i-1}^n \sum_{l=0}^{m+n-k} \sum_{p=0}^{\infty}
\frac{\binom{k}{i}\binom{n}{k}}{\binom{n}{i}}\binom{l+p-1}{p}  \\
& & \quad \times
\binom{m+n-k}{l} \binom{n-i+m-1}{m}(-1)^{l+p+m}q^l(q-1)^p [x]_q^{n-i-m+k+p+l}. \notag
\end{eqnarray}
By (12), we obtain the following theorem.

\begin{theorem}
For $k, n \in \Bbb{Z_+}$ and $i \in \Bbb N$, we have
\begin{eqnarray*}
&&\beta_{i, q}  = \sum_{m=0}^{\infty} \sum_{k=i-1}^n \sum_{l=0}^{m+n-k} \sum_{p=0}^{\infty}
\frac{\binom{k}{i}\binom{n}{k}}{\binom{n}{i}}\binom{l+p-1}{p}  \\
& & \quad \times
\binom{m+n-k}{l} \binom{n-i+m-1}{m}(-1)^{l+p+m}q^l(q-1)^p \beta_{n-i-m+k+p+l, q}. \notag
\end{eqnarray*}
\end{theorem}

\medskip

From (10) and (11),  we note that
\begin{eqnarray}
\frac{\sum_{k=i-1}^{n} \frac{\binom{k}{i}}{\binom{n}{i}}B_{k,n} (x,q)}{([x]_q +[1-x]_q)^{n-i}}
=\sum_{k=0}^{i} q^{\binom{k}{2}} \binom{x}{k}_q [k]_q !S_2 (i, k : q).
\end{eqnarray}
In (3), it is known that
\begin{eqnarray}
\int_{\Bbb Z_p} \binom{x}{n}_q d \mu_q (x) = \frac{(-1)^n}{[n+1]_q}q^{(n+1)-\binom{n+1}{2}}.
\end{eqnarray}
By (13), (14) and Theorem 2, we have
\begin{eqnarray*}
\beta_{n,q}=q \sum_{k=0}^{m} \frac{[k]_q !}{[k+1]_q} (-1)^k S_2 (k, n-k,:q).
\end{eqnarray*}
For $S_2 (n, k :q)$, we see that
\begin{eqnarray}
S_2 (n, k:q)= \frac{1}{(1-q)^k}\sum_{j=0}^{k}(-1)^{k-j}\binom{k+n}{k-j} \binom{j+n}{j}_q, \quad \text{(see [3])},
\end{eqnarray}
and
\begin{eqnarray*}
\binom{n}{k}_q= \sum_{j=0}^{n}\binom{n}{j}(q-1)^{j-k}S_2 (k, j-k :q).
\end{eqnarray*}
By simple calculation, we show that
\begin{eqnarray}
q^{nx}&=& \sum_{k=0}^{n}(q-1)^k q^{\binom{k}{2}} \binom{n}{k}_q [x]_{k, q} \\
&=&\sum_{m=0}^{n}\{\sum_{k=m}^{n}(q-1)^k  \binom{n}{k}_q S_1 (k, m :q)\}[x]_q^m, \notag
\end{eqnarray}
and
\begin{eqnarray}
\int_{\Bbb Z_p} q^{nx} d \mu_q (x) = \sum_{m=0}^{\infty} \binom{n}{m}(q-1)^m
 \beta_{m,q}.\end{eqnarray}
By (16) and (17), we see that
\begin{eqnarray*}
\binom{n}{m}_q= \sum_{k=m}^{n}(q-1)^{-m+k}\binom{n}{k}_q S_1 (k, m :q).
\end{eqnarray*}
Therefore, we obtain the following theorem.

\begin{theorem}
For $k, n \in \Bbb{Z_+}$, we have
\begin{eqnarray*}
B_{k, n}(x,q)=\sum_{m=k}^{n}(q-1)^{-k+m}\binom{n}{m}_q S_1 (m, k :q)[x]_q^k [1-x]_q^{n-k}.
\end{eqnarray*}
\end{theorem}

\medskip

From the definition of the $q$-Stirling numbers of the first kind, we derive
\begin{eqnarray}
q^{\binom{n}{2}}\binom{x}{n}_q [n]_q ! = [x]_{n,q}q^{\binom{n}{2}}=\sum_{k=0}^n S_1 (n, k:q)[x]_q^k.
\end{eqnarray}
By (13) and (18), we obtain the following theorem.

\begin{theorem}
For $k, n \in \Bbb{Z_+}$ and $i \in \Bbb N$, we have
\begin{eqnarray*}
\frac{\sum_{k=i-1}^{n} \frac{\binom{k}{i}}{\binom{n}{i}}B_{k,n} (x,q)}{([1-x]_q +[x]_q)^{n-i}}
=\sum_{k=0}^{i} (\sum_{l=0}^{k}S_1 (n, l : q)  S_2 (i, k : q) [x]_q^l).
\end{eqnarray*}
\end{theorem}

\medskip

By Theorem 2 and Theorem 4, we obtain the following result.
\begin{corollary}
For $i \in \Bbb N$, we have
\begin{eqnarray*}
\beta_{i, q}= \sum_{k=0}^{i}(\sum_{l=0}^{k}S_1(n, l :q )S_2(i, k :q) \beta_{l, q}),
\end{eqnarray*}
where $\beta_{i,q}$ are the $i$-th Carlitz $q$-Bernoulli numbers.
\end{corollary}

\medskip

In [3], the $q$-Bernoulli polynomials of order $k \, (\in \Bbb Z_+)$ are defined by
\begin{eqnarray}
\beta_{n, q}^{(k)}(x)= \frac{1}{(1-q)^n}\sum_{i=0}^{n}(-1)^{i}\binom{n}{i}q^{ix} \underbrace{\int_{\Bbb Z_p}
\cdots \int_{\Bbb Z_p}}_{k-\text{times}}q^{\sum_{l=1}^k (k-l+i)x_l} d \mu_q (x_1) \cdots
d \mu_q (x_k).
\end{eqnarray}
From (19), we note that
\begin{eqnarray*}
\beta_{n, q}^{(k)}= \frac{1}{(1-q)^n}\sum_{i=0}^{n}(-1)^{i}\binom{n}{i}
\frac{(i+k)\cdots(i+1)}{[i+k]_q \cdots [i+1]_q} q^{ix}, \quad (\text{see [3]}).
\end{eqnarray*}
The inverse $q$-Bernoulli polynomial of order $k$ are also defined by
\begin{eqnarray}
\beta_{n, q}^{(-k)}(x)= \frac{1}{(1-q)^n}\sum_{i=0}^{n} \frac{(-1)^{n}\binom{n}{i} q^{xi} }
{\int_{\Bbb Z_p}
\cdots \int_{\Bbb Z_p} q^{\sum_{l=1}^k (k-l+i)x_l} d \mu_q (x_1) \cdots
d \mu_q (x_k)},
\end{eqnarray}
(see [3]). In the special case $x=0$, $\beta_{n, q}^{(k)}=\beta_{n, q}^{(k)}(0)$ are called
the $n$-th $q$-Bernoulli numbers of order $k$ and $\beta_{n, q}^{(-k)}=\beta_{n, q}^{(-k)}(0)$ are called
the $n$-th inverse $q$-Bernoulli numbers of order $k$.

By (20), we see that
\begin{eqnarray}
\beta_{k, q}^{(-n)}&=& \frac{1}{(1-q)^k}\sum_{j=0}^{k}(-1)^{j}\binom{k}{j}
\frac{[j+n]_q \cdots [j+1]_q}{(j+n) \cdots (j+1)} \notag \\
&=& \frac{1}{(1-q)^k}\sum_{j=0}^{k}(-1)^{j} \frac{\binom{k}{j}}{(j+n) \cdots (j+1)}
[n]_q! \binom{j+n}{n}_q \\
&=& \frac{1}{(1-q)^k}\sum_{j=0}^{k}(-1)^{j} \frac{\binom{k+n}{n-j}}{\binom{k+n}{n}}\binom{j+n}{n}_q \frac{[n]_q!}{n!} \notag\\
&=& \frac{[n]_q!}{\binom{k+n}{n} n!} \{ \frac{1}{(1-q)^k} \sum_{j=0}^{k}(-1)^{j}\binom{k+n}{k-j}\binom{j+n}{n}_q \notag \}.
\end{eqnarray}
From (15) and (21), we note that
\begin{eqnarray}
S_2 (n, k :q)=\binom{k+n}{n} \frac{[n]_q !}{n!} \beta_{k,q}^{(-n)}.
\end{eqnarray}
By (13) and (22), we obtain the following theorem.

\begin{theorem}
For $k, n \in \Bbb{Z_+}$ and $i \in \Bbb N$, we have
\begin{eqnarray*}
\frac{\sum_{k=i-1}^{n} \frac{\binom{k}{i}}{\binom{n}{i}}B_{k,n} (x,q)}{([1-x]_q +[x]_q)^{n-i}}
=\sum_{k=0}^{i}q^{\binom{k}{2}} \binom{x}{k}_q [k]_q ! \binom{k+i}{i} \frac{[i]_q !}{i!}
\beta_{k,q}^{(-i)}.
\end{eqnarray*}
\end{theorem}

\medskip

It is not difficult to show that
\begin{eqnarray}
q^{\binom{n}{2}} \binom{x}{n}_q &=& \frac{1}{[n]_q!} \prod_{k=0}^{n-1} ([x]_q - [k]_q) \\
&=& \frac{1}{[n]_q!}  \sum_{k=0}^n (-1)^k [x]_q^{n-k} S_1 (n-1, k :q).
\end{eqnarray}
By Theorem 4 and (23), we obtain the following result.

\begin{corollary}
For $k, n \in \Bbb{Z_+}$ and $i \in \Bbb N$, we have
\begin{eqnarray*}
\frac{\sum_{k=i-1}^{n} \frac{\binom{k}{i}}{\binom{n}{i}}B_{k,n} (x,q)}{([1-x]_q +[x]_q)^{n-i}}
=\sum_{k=0}^{i}\sum_{j=0}^{k}(-1)^j [x]_q^{n-j}S_1 (k-1, j : q)
\binom{k+i}{i} \frac{[i]_q !}{i!}\beta_{k,q}^{(-i)}.
\end{eqnarray*}
\end{corollary}

\end{document}